\documentclass[12pt]{article}
\usepackage{mathrsfs}
\usepackage[justification=centering]{caption}
\usepackage{booktabs}
\usepackage{multirow}
\usepackage{array}
\usepackage[numbers,sort&compress]{natbib}
\usepackage[colorlinks,
            linkcolor=blue,
            anchorcolor=green,
            citecolor=magenta
            ]{hyperref}
\usepackage{mathrsfs,amssymb,amsfonts}
\usepackage{graphicx}
\usepackage{amsmath,amssymb,latexsym,color}
\usepackage[mathscr]{eucal}
\usepackage[numbers]{natbib}
\newtheorem{thm}{Theorem}[section]
\newtheorem{definition}[thm]{Definition}

\newtheorem{lemma}[thm]{Lemma}

\newtheorem{example}{Example}
\newtheorem{remark}{Remark}[section]

\newtheorem{con}{Conjecture}[section]

\newcommand{\proof}{{\it Proof.\quad}}
\newcommand{\qed}{\hfill\Box\medskip}
\usepackage{CJK}
\usepackage{setspace}

\begin{document}

\title{\bf The $g$-extra edge-connectivity of balanced hypercubes}
\author{
Yulong Wei\textsuperscript{a}\footnote{\footnotesize Corresponding author. {\em E-mail address:} weiyulong@tyut.edu.cn (Y. Wei).}
\quad
Rong-hua Li\textsuperscript{b}
\quad
Weihua Yang\textsuperscript{a}
\\
\\{\footnotesize
\textsuperscript{a}\em Department of Mathematics, Taiyuan University of Technology, Taiyuan, 030024, China}
\\{\footnotesize \textsuperscript{b}\em School of Computer Science {\rm \&} Technology, Beijing Institute of Technology, Beijing, 100081, China}}
\date{}
\maketitle

\setlength{\baselineskip}{24pt}

\noindent {\bf Abstract}\quad The $g$-extra edge-connectivity is an important measure for the reliability of interconnection networks. Recently, Yang et al. [Appl. Math. Comput. 320 (2018) 464--473] determined the $3$-extra edge-connectivity of balanced hypercubes $BH_n$ and conjectured that the $g$-extra edge-connectivity of $BH_n$ is $\lambda_g(BH_n)=2(g+1)n-4g+4$ for $2\leq g\leq 2n-1$. In this paper, we confirm their conjecture for $n\geq 6-\dfrac{12}{g+1}$ and $2\leq g\leq 8$, and disprove their conjecture for $n\geq \dfrac{3e_g(BH_n)}{g+1}$ and $9\leq g\leq 2n-1$, where $e_g(BH_n)=\max\{|E(BH_n[U])|\mid U\subseteq V(BH_n), |U|=g+1\}$.


\vskip0.6cm

\section{Introduction}
The topology of interconnection networks can be modeled by a graph $G=(V, E)$ in which a vertex represents a processor and an edge represents a communication link between processors. We refer readers to \cite{BMG,Xu1,Xu2} for terminology and notation unless stated otherwise. Once a network is running, some processors or links might be faulty. An interconnection network without faults is impossible. So the reliability evaluation of interconnection networks is significant.

The traditional {\em edge-connectivity} $\lambda(G)$ is a measurement for the reliability of interconnection networks. However, in real situation, it is a small probability event that all links incident with the same processor fail simultaneously. To overcome this shortcoming, Esfahanian and Hakimi \cite{EH} proposed restricted edge-connectivity. Given a graph $G$, an edge-cut $S\subseteq E(G)$ is called a {\em restricted edge-cut} if there are no isolated vertices in $G-S$. The {\em restricted edge-connectivity} $\lambda'(G)$ is the minimum cardinality of all restricted edge-cuts. Inspired by the restricted edge-connectivity, F$\grave{a}$brega and Foil \cite{FF} proposed the $g$-extra edge-connectivity of a graph. We restate this concept as follows.
\begin{itemize}
\item Given a graph $G$, an edge-cut $F$ is called a {\em $g$-extra edge-cut} if every component of $G-F$ has at least $g+1$ vertices. The {\em $g$-extra edge-connectivity} of $G$, denoted by $\lambda_g(G)$, is the minimum cardinality of all $g$-extra edge-cuts, if exist.
\end{itemize}
A connected graph $G$ is called {\em $\lambda_g$-connected} if $\lambda_g(G)$ exists.

In recent years, the $g$-extra edge-connectivity of a graph has received much attention \cite{L,LX,MS,Y,YFL,YL,ZM,ZXH,ZZ}. For example, Montejano and Sau \cite{MS} proved that given a connected graph $G$ and a positive integer $g$, determining $\lambda_g(G)$ or giving a correct report that $G$ is not $\lambda_g$-connected is NP-hard. Yang \cite{Y} determined that the $1$-extra edge-connectivity of balanced hypercubes $BH_n$ is $\lambda_1(BH_n)=4n-2$ for $n\geq2$. L$\ddot{u}$ \cite{L} showed that $\lambda_2(BH_n)=6n-4$ for $n\geq2$. Li et al. \cite{LX} and Yang et al. \cite{YFL} independently proved that $\lambda_3(BH_n)=8n-8$ for $n\geq2$. In addition, Yang et al.  \cite{YFL} proposed a conjecture about the $g$-extra
edge-connectivity of $BH_n$ as follows.
\begin{con}\label{cj}
Let $BH_n$ be an $n$-dimensional balanced hypercube. Then $\lambda_g(BH_n)=2(g+1)n-4g+4$ for $2\leq g\leq 2n-1$.
\end{con}

Let $e_g(G)=\max\{|E(G[U])|\mid U\subseteq V(G), |U|=g+1\}$, where $G[U]$ is the subgraph of $G$ induced by $U$. In this paper, we confirm their conjecture for $n\geq 6-\dfrac{12}{g+1}$ and $2\leq g\leq 8$, and disprove their conjecture for $n\geq \dfrac{3e_g(BH_n)}{g+1}$ and $9\leq g\leq 2n-1$.

\section{ Balanced hypercubes}\label{2}
In 1997, Wu and Huang proposed balanced hypercubes $BH_n$.
\begin{definition}[\cite{WH}]\label{D0}
An $n$-dimensional balanced hypercube $BH_n=(V(BH_n), E(BH_n))$ has vertex set $V(BH_n)=\{(a_0, a_1, \ldots, a_i, \ldots, a_{n-1})\mid a_i\in\{0, 1, 2, 3\}, 0\leq i\leq n-1\}$. Each vertex $(a_0, a_1, \ldots, a_{i-1}, a_i, a_{i+1}, \ldots, a_{n-1})$ of $BH_n$ has $2n$ neighbors:
\begin{description}
  \item[(1)] $((a_0\pm 1)~{\rm mod~4}, a_1, \ldots, a_{i-1}, a_i, a_{i+1},\ldots, a_{n-1})$,
  \item[(2)] $((a_0\pm 1)~{\rm mod~4}, a_1, \ldots, a_{i-1}, (a_i+(-1)^{a_0})~{\rm mod~4}, a_{i+1},\ldots, a_{n-1})$.
\end{description}
\end{definition}

\begin{figure}[hptb]
  \centering
  \includegraphics[width=8cm]{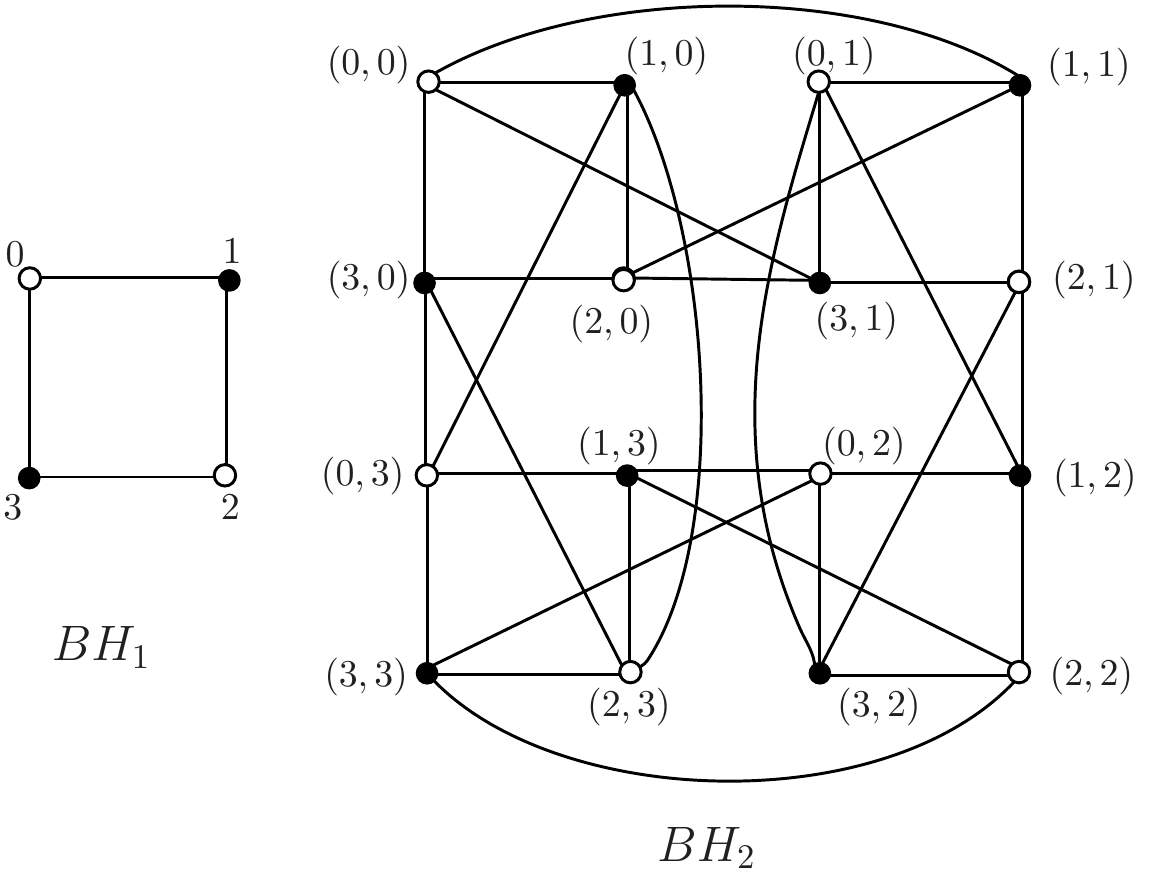}\\
  \caption{Illustration of $BH_1$ and $BH_2$. }\label{DD1}
\end{figure}

Figure \ref{DD1} depicts $BH_1$ and $BH_2$. Clearly, $BH_n$ is a $2n$-regular graph. For a graph $G$ and a vertex $v\in V(G)$, the set of neighbors of $v$ in $G$ is denoted by $N_G(v)$. Some useful properties of $BH_n$ are listed below.

\begin{lemma}[\cite{WH}]\label{pp4}
The balanced hypercube $BH_n$ is bipartite.
\end{lemma}

\begin{lemma}[\cite{Y}]\label{pp2}
Let $u$ be an arbitrary vertex of $BH_n$ for $n\geq1$. Then, for an arbitrary vertex $v$ of $BH_n$, either $|N_{BH_n}(u)\cap N_{BH_n}(v)|=0$, $|N_{BH_n}(u)\cap N_{BH_n}(v)|=2$, or $|N_{BH_n}(u)\cap N_{BH_n}(v)|=2n$. Furthermore, there is exactly one vertex $w$ such that $|N_{BH_n}(u)\cap N_{BH_n}(w)|=2n$.
\end{lemma}

According to Lemma \ref{pp2}, we call the vertex $w$ the {\em equivalent vertex} of $u$, denoted by $u'$,  if $w$ satisfies that $|N_{BH_n}(u)\cap N_{BH_n}(w)|=2n$ in $BH_n$, and $u$ and $u'$ are said to be a pair of equivalent vertices.

The following two lemmas are important observations about the structure of $BH_n$.
\begin{lemma}\label{l1}
The balanced hypercube $BH_n$ is $K_{3,3}$ free.
\end{lemma}
\proof Assume to the contrary that there exists a subgraph $H_1$ of $BH_n$ which is isomorphic to $K_{3,3}$.
By Lemma \ref{pp4}, suppose the bipartite graph $H_1=(X_1, Y_1)$, where $X_1=\{u_1,u_2,u_3\}$ and $Y_1$ are two parts of $H_1$. Since $N_{H_1}(u_1)\cap N_{H_1}(u_2)=Y_1$, $|N_{BH_n}(u_1)\cap N_{BH_n}(u_2)|\geq |N_{H_1}(u_1)\cap N_{H_1}(u_2)|=3>2$. Thus, by Lemma \ref{pp2}, the vertex $u_2$ is the unique equivalent vertex of $u_1$. Similar to the above deduction, we see that the vertex $u_3$ is also the unique equivalent vertex of $u_1$, which contradicts $u_2\neq u_3$.

This completes the proof of Lemma \ref{l1}. $\qed$

Let $\mathcal{F}_g$ be a collection of induced subgraphs of $BH_n$ with $g+1$ vertices and $e_g(BH_n)$ edges for $g\geq2$. By Lemma \ref{pp4}, $H$ is bipartite for any graph $H\in \mathcal{F}_g$.
\begin{lemma}\label{l5}
The vertex set $X$ {\rm(}or $Y${\rm)} must consist of several pairs of equivalent vertices besides at most one vertex for some $H=(X, Y)\in \mathcal{F}_g$.
\end{lemma}
\proof If $|X|=1$, then this lemma holds obviously. Now we consider the case of $|X|\geq2$. Assume to the contrary that there exist two vertices $u, v\in X$ such that their equivalent vertices are not in $X$ for any graph $H=(X, Y)\in \mathcal{F}_g$. Without loss of generality, assume that $|N_{H}(u)|\geq|N_{H}(v)|$. We replace $v$ with $u'$ and obtain an induced subgraph $H'=(X', Y')$ of $BH_n$. If $|N_{H}(u)|>|N_{H}(v)|$, then $|V(H)|=|V(H')|$ and $|E(H)|<|E(H')|$, which contradicts the selection of $H$. If $|N_{H}(u)|=|N_{H}(v)|$ , then $H'\in \mathcal{F}_g$. If $X'$ contains a pair of vertices like $u, v\in X$, then this operation continues until we obtain a graph $H^*=(X^*, Y^*)$ satisfying that $X^*$ consists of several pairs of equivalent vertices besides at most one vertex. Note that $H^*=(X^*, Y^*)\in \mathcal{F}_g$, a contradiction.

By the similar arguments as above, we see that the vertex set $Y$ also consists of several pairs of equivalent vertices besides at most one vertex. $\qed$


The {\em lexicographic product} $G\circ H$ of graphs $G$ and $H$ is defined as the graph with vertex set $V(G)\times V(H)$ and $(u_1, v_1)(u_2, v_2)\in E(G\circ H)$ if and only if $u_1u_2\in E(G)$, or $u_1=u_2$ and $v_1v_2\in E(H)$.
Zhou et al. \cite{ZKF} proved that $BH_n$ is a lexicographic product of a Cayley graph $X_n$ and an empty graph with two vertices. In addition, Zhou et al. \cite{ZWY} showed that $BH_n$ is edge-transitive. Their results are presented as follows.

\begin{lemma}[\cite{ZKF}]\label{lem1}
For each $n\geq1$, $BH_n\cong X_n\circ 2K_1$.
\end{lemma}

\begin{lemma}{\rm\textbf{(}\citep[Page 151]{ZKF}\textbf{)}}\label{lem2}
For $n\geq3$, the girth of $X_n$ is $6$.
\end{lemma}

\begin{lemma}[\cite{ZWY}]\label{main3}
The balanced hypercube is edge-transitive.
\end{lemma}

\section{Main Results}\label{3}
In this section, we will discuss the $g$-extra edge-connectivity of the balanced hypercube $BH_n$ for $2\leq g\leq 2n-1$.

Let $G=(V, E)$ be a graph. For a nonempty proper subset $U\subseteq V$, the set of edges with one end in $U$ and the other end in $\overline{U}=V\setminus U$ is denoted by $[U, \overline{U}]$ and $\partial(U)=|[U, \overline{U}]|$. The {\em $g$-th isoperimetric edge-connectivity} $\gamma_g(G)$ of a graph $G$ was proposed by Hamidoune et al. \cite{HL}. We restate the definition of  $\gamma_g(G)$, that is $\gamma_g(G)=\min\{\partial(U)\mid U\subseteq V, |U|\geq g+1, |\overline{U}|\geq g+1 \}$.
Wang and Li \cite{WL} gave a sufficient condition to ensure a regular edge-transitive graph such that $\lambda_g(G)=\gamma_g(G)$.

\begin{lemma}[\cite{WL}]\label{m1}
Let $G$ be a $k$-regular edge-transitive graph of order $n$ with $k\geq2$, and let $g+1$ be a positive integer. If $n\geq3(g+1)$, then $G$ is $\lambda_g$-connected, and $\lambda_g(G)=\gamma_g(G)$.
\end{lemma}

A graph $G$ satisfying that $\gamma_j(G)=\beta_j(G)$ $(j=0,1,\ldots,g)$ is called $\gamma_g$-{\em optimal}, where $\beta_g(G)=\min\{\partial(U)\mid U\subseteq V, |U|=g+1 \}$.
Zhang \cite{ZZ} gave a sufficient condition for a regular edge-transitive graph to be $\gamma_g$-optimal.

\begin{lemma}[\cite{ZZ}]\label{m2}
Let $g+1$ be a positive integer, and $G$ a connected $k$-regular edge-transitive graph with $k\geq\dfrac{6e_g(G)}{g+1}$. Then $G$ is $\gamma_g$-optimal.
\end{lemma}

The following lemma gives a lower bound of $e_g(BH_n)$ for $2\leq g\leq 2n-1$.

\begin{figure}[hptb]
  \centering
  \includegraphics[width=8cm]{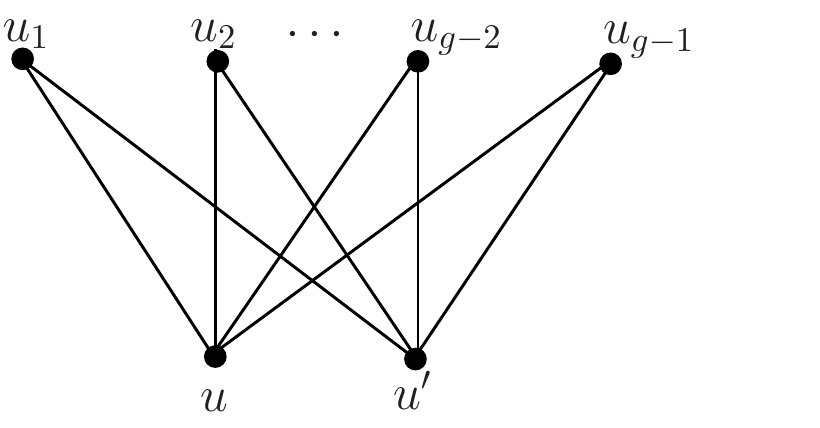}\\
  \caption{Illustration of $BH_n[A_g]$. }\label{f2}
\end{figure}

\begin{lemma}\label{l2}
The balanced hypercube $BH_n$ satisfies that $e_g(BH_n)\geq2g-2$ for $2\leq g\leq 2n-1$.
\end{lemma}
\proof Suppose that $u=(0,0,\ldots,0)$, $u'=(2,0,\ldots,0)$, $u_1=(1,0,\ldots,0)$, $u_2=(3,0,\ldots,0)$, $u_{2i-1}=(1,\overbrace{0,\ldots,0}^{i-2},1,0,\ldots,0)$ and $u_{2i}=(3,\overbrace{0,\ldots,0}^{i-2},1,0,\ldots,0)$ for $2\leq i\leq g-1$ are some vertices of $BH_n$. Let $A_g=\{u, u'\}\cup\{u_i\mid 1\leq i\leq g-1\}$. By Definition \ref{D0}, we know that the induced subgraph $BH_n[A_g]$ is isomorphic to $K_{2, g-1}$ (see Figure \ref{f2}). Therefore, $e_g(BH_n)\geq2g-2$ for $2\leq g\leq 2n-1$.  $\qed$

Now, we determine $e_g(BH_n)$ for $2\leq g\leq 8$.

\begin{lemma}\label{l3}
The balanced hypercube $BH_n$ satisfies that $e_g(BH_n)=2g-2$ for $2\leq g\leq 8$.
\end{lemma}
\proof
By Lemma \ref{l2}, we only need to prove that  $e_g(BH_n)\leq 2g-2$ for $2\leq g\leq8$. Let $H$ be an induced subgraph of $BH_n$ with $|V(H)|=g+1$ and $|E(H)|=e_g(BH_n)$. By Lemma \ref{pp4}, suppose the bipartite graph $H=(X, Y)$, where $X$ and $Y$ are two parts of $H$. We divide our discussion into five cases.

{\em Case 1.} $2\leq g\leq4$.

In this case, $|X|+|Y|=|V(H)|=g+1$. Hence, $$|E(H)|\leq|X|\cdot|Y|\leq \left\lfloor\dfrac{g+1}{2}\right\rfloor\cdot\left\lceil\dfrac{g+1}{2}\right\rceil=2g-2. $$

{\em Case 2.} $g=5$.

In this case, $|X|+|Y|=|V(H)|=6$. Hence, $|E(H)|\leq|X|\cdot|Y|\leq 3\times3=9$. If $|E(H)|=9$, then $H$ is isomorphic to $K_{3,3}$, which contradicts Lemma \ref{l1}. Thus, $|E(H)|\leq8=2g-2$.

{\em Case 3.} $g=6$.

In this case, $|X|+|Y|=|V(H)|=7$. Hence, $|E(H)|\leq|X|\cdot|Y|\leq 3\times4=12$. If $|E(H)|=12$, then $H$ is isomorphic to $K_{3,4}$, which contradicts Lemma \ref{l1}. If $|E(H)|=11$, then $H$ is isomorphic to $K_{3,4}-e$ for some $e\in E(K_{3,4})$, which also contradicts Lemma \ref{l1}. Thus, $|E(H)|\leq10=2g-2$.

{\em Case 4.} $g=7$.

In this case, $|X|+|Y|=|V(H)|=8$. Hence, $|E(H)|\leq|X|\cdot|Y|\leq 4\times4=16$. Note that $|E(K_{2,6})|=12=2g-2$. Let $E_0\subseteq E(K_{3,5})$ and $E_1\subseteq E(K_{4,4})$. If $H=K_{3,5}-E_0$ and $|E_0|\leq2$, then $H$ contains a subgraph isomorphic to $K_{3, 3}$, which contradicts Lemma \ref{l1}. If $H=K_{4,4}-E_1$ and $|E_1|\leq3$, then by Lemma \ref{l5}, suppose that $X=\{u_1, u_1', u_2, u_2'\}$ and $Y=\{v_1, v_1', v_2, v_2'\}$. Since $H=(X, Y)$ is an induced subgraph of $BH_n$, $|E_1|=0$. Then $H$ is isomorphic to $K_{4, 4}$, which contradicts Lemma \ref{l1}. Thus, $|E(H)|\leq12=2g-2$.

{\em Case 5.} $g=8$.

In this case, $|X|+|Y|=|V(H)|=9$. Hence, $|E(H)|\leq|X|\cdot|Y|\leq 4\times5=20$. Note that $|E(K_{2,7})|=14=2g-2$. Let $E_0\subseteq E(K_{3,6})$ and $E_1\subseteq E(K_{4,5})$.  Now, we only need to discuss the following two cases.

{\em Case 5.1.} $H=K_{3,6}-E_0$ and $|E_0|\leq3$.

Note that $K_{3,6}-E_0$ with $|E_0|\leq3$ contains a subgraph isomorphic to $K_{3,3}$, which contradicts Lemma \ref{l1}.

{\em Case 5.2.} $H=K_{4,5}-E_1$ and $|E_1|\leq5$.

By Lemma \ref{l5}, without loss of generality, suppose that $X=\{u_1, u_1', u_2, u_2'\}$ and $Y=\{v_1, v_1', v_2, v_2', v\}$. Note that $H=(X, Y)$ is an induced subgraph of $BH_n$.  Then $|E_1|$ can not be an odd integer. If $|E_1|\in\{0, 2\}$, then $K_{4,5}-E_1$ contains a subgraph isomorphic to $K_{3,3}$, which contradicts Lemmas \ref{l1}. If $|E_1|=4$ and $u_iv_j\in E_1$ for some $i, j\in\{1, 2\}$, then edges $u_i'v_j, u_iv_j', u_i'v_j'\in E_1$. Thus, $E_1=\{u_iv_j, u_i'v_j, u_iv_j', u_i'v_j'\}$ and $H=K_{4,5}-E_1$ contains a subgraph isomorphic to $K_{3,3}$, which contradicts Lemmas \ref{l1}. If $|E_1|=4$ and $u_iv_j\notin E_1$ for all $i, j\in\{1, 2\}$, then $E_1=\{u_1v, u_2v, u_1'v, u_2'v\}$. Therefore, $H=K_{4,5}-E_1$ contains a subgraph isomorphic to $K_{3,3}$, which contradicts Lemmas \ref{l1}.

Thus, $|E(H)|\leq14=2g-2$.

So, $e_g(BH_n)=2g-2$ for $2\leq g\leq8$. $\qed$

The following lemma gives a lower bound of $e_g(BH_n)$ for $9\leq g\leq 2n-1$, which will be used to disprove Conjecture \ref{cj} for $9\leq g\leq 2n-1$.

\begin{lemma}\label{l8}
The balanced hypercube $BH_n$ satisfies that $e_g(BH_n)>2g-2$ for $9\leq g\leq 2n-1$.
\end{lemma}
\proof To prove this lemma, it suffices to construct a subgraph of $BH_n$ with $g+1$ vertices and at least $2g-1$ edges.

By Lemma \ref{lem2}, the girth of $X_n$ is $6$ for $n\geq3$. Suppose that $\overline{C_6}$ is a cycle of $X_n$ with six vertices. Let $H_0=\overline{C_6}\circ 2K_1$ be a subgraph of $BH_n$. Since $BH_n$ is connected, by Lemma \ref{lem1}, $X_n$ is a connected graph for $n\geq1$. Let $\overline{U_t}$ be a connected subgraph of $X_n$ with $|V(\overline{U_t})|=t\geq6$ satisfying that $\overline{U_t}$ is a unicyclic graph which contains $\overline{C_6}$. Then $|E(\overline{U_t})|=t\geq6$. Now, we distinguish the following four cases.

{\em Case 1.} $g=9$.

We consider the graph $H_0-\{u, v\}$, where $u, v\in V(H_0)$, $u'\neq v$ and $uv\in E(H_0)$. Note that
$|V(H_0-\{u, v\})|=g+1$ and $|E(H_0-\{u, v\})|=2g-1$. Then $e_g(BH_n)\geq|E(H_0-\{u, v\})|>2g-2$.

{\em Case 2.} $g=10$.

We consider the graph $H_0-v$, where $v\in V(H_0)$. Note that
$|V(H_0-v)|=g+1$ and $|E(H_0-v)|=2g$. Then $e_g(BH_n)\geq|E(H_0-v)|>2g-2$.

{\em Case 3.} $g$ is an odd integer with $g\geq11$.

Since $g$ is an odd integer with $g\geq11$, we have $\dfrac{g+1}{2}\geq6$. We consider the graph $H_g=\overline{U_{\frac{g+1}{2}}}\circ 2K_1$ as a subgraph of $BH_n$. Note that $|V(H_g)|=g+1$ and $|E(H_g)|=2g+2$. Then $e_g(BH_n)\geq|E(H_g)|>2g-2$.

{\em Case 4.} $g$ is an even integer with $g\geq12$.

Since $g$ is an even integer with $g\geq12$, $g-1$ is an odd integer with $g-1\geq11$. We consider the graph $H_{g-1}=\overline{U_{\frac{g}{2}}}\circ 2K_1$ as a subgraph of $BH_n$. Pick a vertex $u$ from $V(BH_n)\setminus V(H_{g-1})$. Note that $|V(H_{g-1})\cup \{u\}|=g+1$ and $|E(BH_n[V(H_{g-1})\cup \{u\}])|\geq|E(H_{g-1})|=2g$. Then $e_g(BH_n)\geq|E(BH_n[V(H_{g-1})\cup \{u\}])|>2g-2$.

As mentioned above, we obtain the desired result. $\qed$

Now, we give the proof of our main theorem.

\begin{thm}\label{l4}
The $g$-extra edge-connectivity of balanced hypercubes $BH_n$ is $\lambda_g(BH_n)=2(g+1)n-4g+4$ for $n\geq 6-\dfrac{12}{g+1}$ and $2\leq g\leq 8$. In addition, $\lambda_g(BH_n)<2(g+1)n-4g+4$ for $n\geq \dfrac{3e_g(BH_n)}{g+1}$ and $9\leq g\leq 2n-1$.
\end{thm}
\proof By Lemma \ref{main3}, $BH_n$ is edge-transitive. Note that $|V(BH_n)|=2^{2n}\geq6n\geq3(g+1)$ for $n\geq2$. By Lemma \ref{m1}, $\lambda_g(BH_n)=\gamma_g(BH_n)$ for $2\leq g\leq 2n-1$.

By Lemma \ref{l3}, $e_g(BH_n)=2g-2$ for $2\leq g\leq 8$. Since $n\geq 6-\dfrac{12}{g+1}$, we have $2n\geq \dfrac{6(2g-2)}{g+1}=\dfrac{6e_g(BH_n)}{g+1}$ for $2\leq g\leq 8$. By Lemma \ref{m2}, $BH_n$ is $\gamma_g$-optimal. Thus, $\gamma_g(BH_n)=\beta_g(BH_n)$ for $n\geq 6-\dfrac{12}{g+1}$ and $2\leq g\leq 8$.
Since $\beta_g(BH_n)=2n(g+1)-2e_g(BH_n)=2n(g+1)-2(2g-2)=2(g+1)n-4g+4$, we have $\lambda_g(BH_n)=2(g+1)n-4g+4$ for $n\geq 6-\dfrac{12}{g+1}$ and $2\leq g\leq 8$.

By Lemma \ref{l8}, $e_g(BH_n)>2g-2$ for $9\leq g\leq 2n-1$. By Lemma \ref{m2}, we have  $\gamma_g(BH_n)=\beta_g(BH_n)=2n(g+1)-2e_g(BH_n)<2n(g+1)-2(2g-2)=2(g+1)n-4g+4$ for $2n\geq \dfrac{6e_g(BH_n)}{g+1}$ and $9\leq g\leq 2n-1$. Therefore, $\lambda_g(BH_n)<2(g+1)n-4g+4$ for $n\geq \dfrac{3e_g(BH_n)}{g+1}$ and $9\leq g\leq 2n-1$. $\qed$

\section{Conclusions}\label{4}
The $g$-extra edge-connectivity is an important measure for the reliability of interconnection networks. We establish the $g$-extra edge-connectivity of balanced hypercubes $BH_n$, that is $\lambda_g(BH_n)=2(g+1)n-4g+4$ for $n\geq 6-\dfrac{12}{g+1}$ and $2\leq g\leq 8$, which partially confirms Conjecture \ref{cj}. This result can provide a more accurate measurement of edge fault tolerance of balanced hypercubes. Meanwhile, we prove that $\lambda_g(BH_n)<2(g+1)n-4g+4$ for $n\geq \dfrac{3e_g(BH_n)}{g+1}$ and $9\leq g\leq 2n-1$, which disproves Conjecture \ref{cj} for any $n$ and $g$ with $n\geq \dfrac{3e_g(BH_n)}{g+1}$ and $9\leq g\leq 2n-1$.

\section*{Acknowledgement}
Y. Wei's research is supported by the Natural Science Foundation of Shanxi Province (No. 201901D211106). W. Yang's research is supported by the National Natural Science Foundation of China (No. 11671296).

\noindent

\end{document}